\documentclass[12pt]{article}
\usepackage{amsfonts, amsthm, amsmath}

\newtheorem{theorem}{Theorem}
\newtheorem{lemma}[theorem]{Lemma}

\def\FF{\mathbb{F}}
\def\PP{\mathbb{P}}

\def\Abar{\overline{A}}

\def\Qbar{\overline{Q}}

\def\Tr{\mathrm{Tr}}
\def\ZZ{\mathbb{Z}}
\def\idm{\mathfrak{m}}

\begin{document}

\title{Counting Points on
Hyperelliptic Curves using Monsky-Washnitzer Cohomology}
\author{Kiran S. Kedlaya}
\date{November 20, 2001}
\maketitle

\begin{abstract}
We describe an algorithm for counting points on an arbitrary hyperelliptic
curve over
a finite field $\FF_{p^n}$ of odd characteristic,
using Monsky-Washnitzer cohomology to
compute a $p$-adic approximation to the characteristic polynomial of
Frobenius. For fixed $p$,
the asymptotic running time for a curve of genus $g$ over $\FF_{p^n}$
with a rational Weierstrass point
is $O(g^{4+\epsilon} n^{3+\epsilon})$.
\end{abstract}

\section{Introduction}

An important problem in computational algebraic geometry is the enumeration
of points on algebraic varieties over finite fields, or more precisely
the determination of their zeta functions.
Much work so far on this problem has focused on curves of genus 1.
Initial approaches, like the Shanks-Mestre method \cite[Section~7.4.3]{cohen},
yield algorithms with exponential running time in the length of the input
data (which is roughly the logarithm of the field size).
Schoof \cite{schoof} gave an algorithm for counting points on a genus 1 curve
over $\FF_q$ which is polynomial in $\log(q)$; this algorithm was improved
by Atkin and Elkies. For fields of fixed (or at least small)
characteristic, an algorithm
given by Satoh \cite{satoh} has smaller asymptotic running time than
Schoof's algorithm; an implementation is described in detail in \cite{fgh}.

Extending the aforementioned methods to curves of higher genus
has to date yielded unsatisfactory results. The Shanks-Mestre method is
exponential both in the field size and in the genus.
Schoof's algorithm, which is roughly to compute the characteristic
polynomial of Frobenius modulo many small primes, can be generalized in
principle to higher genus, as noted by Pila \cite{pila}. However,
using the method in practice requires producing
explicit equations for the Jacobian of the curve,
which is already nontrivial in genus 2 and probably hopeless in general.
Satoh's method,
which is to compute the Serre-Tate canonical lift, runs into a similar
obstruction: the Serre-Tate lift of a Jacobian need not itself be a Jacobian,
so computing with it is difficult. Satoh has proposed working instead with
the formal group of the Jacobian. This is possible in principle, as the
formal group can be expressed in terms of data on the curve, but the result
again seems to be exponential in the genus.

In this paper, we develop an
algorithm for counting points on hyperelliptic curves over finite
fields of odd characteristic, which is
polynomial in the genus of the curve.
Our approach is to compute in the Monsky-Washnitzer (dagger)
cohomology of an affine
curve, which is essentially the de~Rham cohomology of a lift of the curve to
characteristic zero, endowed with an action of Frobenius. The action
of Frobenius can be $p$-adically approximated efficiently using certain
power series.

As the approach is $p$-adic, the method shares with Satoh's algorithm
the nature of its dependence on the input parameters. Namely,
both algorithms are polynomial
in the degree of the finite field over the prime
field (with the same exponent), but are polynomial in the order
of the prime field rather than in its logarithm. Additionally, our 
algorithm is polynomial in the genus of the hyperelliptic curve.
To be specific, the running time of the algorithm is on the order
of $g^{4+\epsilon} n^{3+\epsilon}$, where $n$ is the field degree and $g$
the genus, assuming that the curve has a rational Weierstrass point.
(One should be able to achieve the same running time even without a
rational Weierstrass point, but we have not checked this.)

The strategy of counting points on a variety by computing in de~Rham
cohomology on a lift seems to be quite broadly applicable. In particular,
there is no reason why it could not be applied to more general curves,
or even to higher dimensional varieties (e.g.,
hypersurfaces in toric varieties). In fact, a related method has been
introduced by Lauder and Wan \cite{lw}, who use Dwork's
trace formula to give a $p$-adic
algorithm for computing the zeta function of an
arbitrary variety over a finite field. It is unclear how practical it will be
to implement that algorithm; Lauder and Wan themselves suggest reinterpreting
it in terms of a $p$-adic cohomology theory to make it easier to implement.

\section{Overview of $p$-adic cohomology}
We briefly recall the formalism of Monsky-Washnitzer cohomology,
as introduced by Monsky and Washnitzer \cite{mw}, \cite{m2}, \cite{m3}
and refined by van der Put \cite{vdp}; details omitted here can be
found therein.
We first set some notations.
Let $k$ be a perfect field of characteristic $p>0$
(which for us will always be a finite field), $R$
a complete mixed characteristic discrete valuation ring with residue field
$k$ (e.g., the ring of Witt vectors $W(k)$), and $\idm$ the maximal ideal of $R$.
Let $K$ be the fraction field of $R$.

Let $X$ be a smooth affine variety over $k$,
$\Abar$ the coordinate ring of $X$, and $A$ a smooth $R$-algebra
with $A \otimes_R k \cong \Abar$.
Ordinarily, $A$ will not admit a lift of the
absolute Frobenius morphism on $\Abar$, but its $\idm$-adic
completion $A^\infty$
will. Working with $A^\infty$ is not satisfactory, however,
because the de~Rham cohomology of $A^\infty$
is larger than that of $A$. The trouble is that
the limit of exact differentials
need not be exact: for example, if $A = R[x]$, then
the sum $\sum_{n=0}^\infty p^n x^{p^n-1}\,dx$ defines a differential over $A$
which is the limit of exact differentials but is not itself exact.

To remedy the situation, Monsky and Washnitzer work with a subring of
$A^\infty$, consisting of series which converge fast enough
that their integrals also converge. Namely, fix $x_1, \dots, x_m \in
A^\infty$ whose images generate $\Abar$ over $k$. Monsky and Washnitzer
define the \emph{weak completion}
$A^\dagger$ of $A$ as the subring of $A^\infty$ consisting of elements $z$
representable, for some real number $c$,
as $\sum_{n=0}^\infty a_n P_n(x_1, \dots, x_n)$, with $a_n \in \idm^n$ and
$P_n$ an $n$-variate polynomial of total degree
at most $c(n+1)$. One can show that the weak completion depends, up
to noncanonical isomorphism, only on $\Abar$.

Monsky and Washnitzer then define the dagger cohomology groups
$H^i(\overline{A}; K)$ as the cohomology groups of the de~Rham
complex over $A^\dagger \otimes_R K$. Namely, let $\Omega$
denote the $A^\dagger$-module
of differential forms over $K$, generated by symbols
$dx$ for $x \in A^\dagger \otimes_R K$ and subject to the relations
$d(xy) = x\,dy + y\,dx$ for all $x$ and $y$, and $dx =0$ for $x \in K$.
Then the map $d: \wedge^i \Omega \to \wedge^{i+1} \Omega$ given
by
\[
x\,dy_1 \wedge \cdots \wedge dy_i = dx \wedge dy_1 \wedge \cdots
\]
satisfies $d \circ d = 0$, and thus makes the $\Omega^i$ into a complex,
whose cohomology at $\Omega^i$ we call $H^i(\overline{A}; K)$;
this group is in fact a $K$-vector space.
This construction is clearly functorial with respect to maps on dagger
rings; in particular, if $\phi$ is an endomorphism of $A^\dagger$, it
induces an endomorphism $\phi_*$ on the cohomology groups.

The point of this construction is that these cohomology groups satisfy
the following Lefschetz fixed point formula.
See van~der~Put \cite[4.1]{vdp} for a proof.
\begin{theorem}[Lefschetz fixed point formula] \label{lefschetz}
Let $\overline{A}$ be smooth and integral of dimension $n$
over $\FF_q$. Suppose the weak completion $A^\dagger$ of a lift of
$\overline{A}$ admits an endomorphism $F$ lifting the $q$-power Frobenius on
$\overline{A}$.
Then the number
of homomorphisms $\overline{A} \to \FF_q$ equals
\[
\sum_{i=0}^n (-1)^i \Tr(q^n F_*^{-1} | H^i(\overline{A}; K)).
\]
\end{theorem}

In the original work of Monsky-Washnitzer, it was unknown whether the
cohomology groups were necessarily finite dimensional as vector
spaces over $K$; thus in the fixed
point formula, the fact that the operator $F_*^{-1}$ has a trace is
a nontrivial part of the result. It was later shown by Berthelot
\cite{berthelot} that the vector spaces are indeed finite dimensional.
Thus we can compute the traces in the fixed point formula by working
in finite dimensional vector spaces.

Summing up, we have the following general strategy for
computing the zeta function
of a smooth projective variety $X$ over a finite field $\FF_q$;
we flesh out this strategy in the particular case at hand in the
rest of the paper.
Choose an affine subvariety $U$ of $X$, then compute the zeta function of
$X-U$, which is a closed subvariety of $X$ of lower dimension. Then
compute the action of a lift of Frobenius on the cohomology groups of
$U$; since one cannot exactly represent all elements of $W(\FF_q)$, the
action can only be computed to a certain $p$-adic precision. The
net result is a $p$-adic approximation of the zeta function; by using
enough precision, one can get a good enough approximation that the
Riemann hypothesis component of the Weil conjectures uniquely
determines the zeta function from this approximation.

\section{Cohomology of hyperelliptic curves}


In this section, let $p$ be an odd prime. 
We describe the Monsky-Washnitzer
cohomology of a hyperelliptic curve over a field of characteristic
$p$ in a concrete manner, suitable for explicit computation of its
zeta function; we will 
explicitly describe the computation in the next section.

We begin by setting notation for this section and the next.
Let $\Qbar(x)$ be a polynomial of degree $2g+1$
over $\FF_{q}$ without repeated roots, so that
the closure in the projective plane of the affine curve
$y^2 = \Qbar(x)$ is a smooth hyperelliptic curve $C$
of genus $g$ with a rational Weierstrass point. (One can handle the case
where there is no rational Weierstrass point by similar methods, but we
omit the details here.)
Let $C'$ be the affine curve obtained from $C$ by deleting the support
of the divisor of $y$ (that is, the point at infinity and the Weierstrass
points); then the coordinate ring $\Abar$ of $C'$ is 
$\FF_{q}[x,y,y^{-1}]/(y^2 - \Qbar(x))$.
Let $A = W(\FF_{q})[x,y,y^{-1}]/(y^2-Q(x))$ and
let $A^\dagger$ be the weak completion of $A$.

Before proceeding further, we give an explicit description of $A^\dagger$.
Namely, let $v_p$ denote the $p$-adic valuation on $W(\FF_q)$,
and extend this norm to polynomials as follows:
if $P(x) = \sum a_i x^i$, define $v_p(P) = \min_i \{v_p(a_i)\}$.
Then the elements of $A^\dagger$ can be viewed as
series $\sum_{n=-\infty}^\infty
(S_n(x) + T_n(x)y)y^{2n}$, where $S_n$ and $T_n$ are polynomials of degree
at most $2g$, such that 
\[
\liminf_{n \to \infty} \frac{v_p(S_n)}{n}, \quad
\liminf_{n \to \infty} \frac{v_p(T_n)}{n}, \quad
\liminf_{n \to \infty} \frac{v_p(S_{-n})}{n}, \quad
\liminf_{n \to \infty} \frac{v_p(T_{-n})}{n}
\]
are all positive.

We can lift the $p$-power Frobenius to
an endomorphism $\sigma$ of $A^\dagger$
by defining it as the canonical Witt vector Frobenius on $W(\FF_q)$, 
then extending to
$W(\FF_q)[x]$ by mapping $x$ to $x^p$,
and finally setting
\begin{align*}
y^\sigma &= y^p \left( 1 + \frac{Q(x)^\sigma - Q(x)^p}{Q(x)^p} \right)^{1/2} \\
&= y^p \sum_{i=0}^\infty \binom{1/2}{i} \frac{(Q(x)^\sigma - Q(x)^p)^i}{Q(x)^{pi}} \\
&= y^p \sum_{i=0}^\infty \frac{(1/2)(1/2 - 1) \cdots (1/2-i+1)}{i!}
\frac{(Q(x)^\sigma - Q(x)^p)^i}{Q(x)^{pi}}
\end{align*}
and $(y^{-1})^\sigma = (y^\sigma)^{-1}$.
Let $F = \sigma^{\log_p q}$; then $F$ is a lift of the $q$-power Frobenius,
so we may apply the Lefschetz fixed point formula to it and use the
result to compute the zeta function of $C$. We now describe how this is done.

The de~Rham cohomology of $A$ splits into eigenspaces under the
hyperelliptic involution: a positive eigenspace generated by
$x^i\,dx / y^2$ for $i=0, \dots, 2g-1$, and a negative eigenspace
generated by $x^i\,dx/y$ for $i=0, \dots, 2g-1$. Indeed, any form can be
written as $\sum_{n=-\infty}^\infty
 \sum_{i=0}^{2g-1} a_{i,n} x^i\,dx/y^n$, and the relation
\[
\frac{B'(x)\,dx}{y^{s}} \equiv \frac{s B(x)\,dy}{y^{s+1}}
= \frac{s B(x)Q'(x)\,dx}{2y^{s+2}}
\]
(which follows from the equality $2y\,dy = Q'(x)\,dx$)
can be used to consolidate everything into the $n=1$ and $n=2$ terms.
Specifically, when $s >1$, we can write an arbitrary polynomial $B(x)$
as $R(x)Q(x) + S(x) Q'(x)$ for some polynomials $R, S$ (since $Q$ has
no repeated roots), and then write
\[
\frac{B(x)\,dx}{y^s} \equiv \frac{R(x)\,dx}{y^{s-2}} + \frac{2S'(x)\,dx}{(s-2)y^{s-2}}.
\]
On the other side, a differential $B(x)\,dx/y$ with $B$ a polynomial of 
degree greater than $2g$ can be reduced using the identity
$[S(x) Q'(x) + 2 S'(x) Q(x)]dx/y \equiv 0$.
For $S(x) = x^{m-2g}$, the expression in brackets has degree $m$
and leading term $(2g+1) + 2(m-2g) = 2m-2g+1 \neq 0$, so a suitable multiple
can be subtracted from $B$ to reduce its degree.

To carry out provably
correct computations, we need explicit estimates on the denominators 
introduced by the aforementioned reduction process.
We now prove a lemma that provides the needed estimate. (The approach is
similar to that of the proof of \cite[Lemma~4.1]{m2}.)

\begin{lemma} \label{lem:est}
Let $A(x)$ be a polynomial over $W(\FF_q)$ of degree at most $2g$. Then 
for $m \geq 0$, the
reduction of $\omega = A(x)\,dx/y^{2m+1}$ becomes integral upon
multiplication by $p^{\lfloor \log_p (2m+1) \rfloor}$.
\end{lemma}
\begin{proof}
Let $B(x)\,dx/y$ be the reduction of $A(x)\,dx/y^{2m+1}$, and
$f$ the function such that $df = A(x)\,dx/y^{2m+1} - B(x)\,dx/y$.
Write $f = \sum_{j=0}^{m} F_j(x)/y^{2j+1}$ where each $F_j$ has
degree at most $2g$.
Let $r_0, \dots, r_{2g}$ be the roots of $Q(x)$ over $W(\overline{\FF_q})$
and $T_0, \dots, T_{2g}$
the corresponding points on the curve $y^2 = Q(x)$. Then $f$ has poles
at $T_0, \dots, T_{2g}$ and possibly at infinity.
 
Let $R_i$ be the completion of the local ring of $W(\overline{\FF_q})[x,y] /
(y^2-Q(x))$ at
$T_i$, and let $K_i$ the fraction field of $R_i$; then the maximal ideal of
$R_i$ is generated by $y$, and within $R_i$,
$x$ can be written as a power series in $y$ with integral coefficients.
Then the image of $df$ in the module $\Omega_{K_i/W(\overline{\FF_q})}$
of differentials
can be written as
$\sum_{k=-m}^\infty a_{ik} y^{2k-2}\,dy$, and the $a_{ik}$ are integral
for $k<0$ (since they coincide with the corresponding coefficients in 
the expansion of $\omega$).

The map $d$ commutes with the passage to the completed local ring, 
so the image of $f$ in $K_i$ is equal to $\sum_{k=-m}^\infty
a_{ik} y^{2k-1}/(2k-1)$. Now note that $f - \sum_{k=-m}^{-j-1}
F_{-k}(x) y^{2k-1}$ has a pole of order at most $2j+1$ at each $T_i$,
and its image in $K_i$ has leading term $F_{j}(r_i) y^{-2j-1}$.
Consequently, if $n$ is an integer such that $n a_{ik}/(2k-1)$ is integral
for $i=0,\dots,2g$ and $k=-1,\dots,-m$, then $nf$ is integral. Specifically,
we have that $nF_{-m}(r_i)$ is integral for $i=1, \dots, 2g+1$;
since the $r_i$ are distinct modulo $p$, that implies that $nF_{-m}(x)$ is
integral. Applying the same argument to $nf - nF_{-m}(x)$, we deduce that
$nF_{-m+1}(x)$ is integral, and so forth.

In particular, we may take $n = p^{\lfloor \log_p (2m+1) \rfloor}$.
Then $nf$ is integral, as is the reduction of $n \omega$,
which yields the desired conclusion.
\end{proof}
One can make the following analogous assertion for the reduction process
in the other direction, using the local ring at infinity instead of at
the $T_i$. We omit the proof.
\begin{lemma}
Let $A(x)$ be a polynomial over $W(\FF_q)$ of degree at most $2g$. Then 
for $m \geq 0$, the
reduction of $\omega = A(x)y^{2m+1}\,dx$ becomes integral upon
multiplication by $p^{\lfloor \log_p (2m+1) \rfloor}$.  
\end{lemma}
In particular, the basis we have chosen of the de Rham cohomology of $A$
is also a basis of $H^1(\Abar; K)$.

Let $F = \sigma^{\log_p q}$ denote the $q$-power Frobenius.
By the Lefschetz fixed point formula (Theorem~\ref{lefschetz})
applied to $C'$ and its image $\PP^1$ under quotienting by the
hyperelliptic involution, we have
\begin{align*}
\#C(\FF_{q^i}) - 2g &= \#C'(\FF_{q^i})\\
&= \Tr(q^iF^{-i}, H^0(\Abar; K)) - \Tr(q^iF^{-i}, H^1(\Abar; K)) \\
&= q^i - \Tr(q^iF^{-i}, H^1(\Abar; K)_{+}) -
\Tr(q^iF^{-i}, H^1(\Abar; K)_{-}) \\
&= \Tr(q^iF^{-i}, H^0(\Abar; K)_{+}) -  \Tr(q^iF^{-i}, H^1(\Abar; K)_{+}) - \Tr(q^iF^{-i}, H^1(\Abar; K)_{-}) \\
&= \#\PP^1(\FF_{q^i}) - \Tr(q^iF^{-i}, H^1(\Abar; K)_{-}) \\
&= (q^i+1-2g) - \Tr(q^iF^{-i}, H^1(\Abar; K)_{-}).
\end{align*}
Thus $q+1-\#C(\FF_{q^i})$ equals the trace of $q^iF^{-i}$ on the
negative eigenspace of $H^1(\Abar; K)$.

By the Weil conjectures (see \cite[Appendix~C]{hartshorne} for details),
there exists a polynomial
$x^{2g} + a_1 x^{2g-1} + \cdots + a_{2g}$ whose roots
$\alpha_1, \dots, \alpha_{2g}$ satisfy
$\alpha_j \alpha_{g+j} = q$ for $j=1, \dots, g$,
$|\alpha_j| = \sqrt{q}$ for $j=1, \dots, 2g$, and 
\[
q+1-\#C(\FF_{q^i}) = \sum_{j=1}^{2g} \alpha_j^i
\]
for all $i>0$.
Thus the eigenvalues of $qF^{-1}$ on $H^1(\Abar; K)_{-}$ are precisely the
$\alpha_i$, as are the eigenvalues of $F$ itself. 
Since $a_{i} = a_{2g-i}$, it suffices to determine $a_1, \dots, a_g$.
Moreover,
$a_i$ is the sum of $\binom{2g}{i}$ $i$-fold products of eigenvalues of
Frobenius, so for $i=1, \dots, g$,
\[
|a_i| \leq \binom{2g}{i} q^{i/2} \leq 2^{2g} q^{g/2}.
\]
Thus to determine the zeta function,
it suffices to compute the action of $F$ on a suitable
basis of $H^1(\Abar; K)_{-}$ modulo
$p^{N_1}$ for $N_1 \geq (g/2)n + (2g+1) \log_p 2$.
Thanks to Lemma~\ref{lem:est},
we can determine explicitly how much computation is needed to 
determine this action.

The action of the $p$-power Frobenius $\sigma$ on differentials is given by 
\[
\left( \frac{A(x)\,dx}{y^{2k+1}} \right)^\sigma
 = \frac{p A(x)^\sigma x^{p-1}\,dx}{y^{p(2k+1)}}
\left( 1 + \frac{pE(x)}{y^{2p}} \right)^{-(2k+1)/2},
\]
where we set $pE(x) = Q(x)^\sigma - Q(x)^p$. We can rewrite this expression
as a power series $\sum_i A_i(x) y^{-2i-1}\,dx$, where each
polynomial $A_i(x)$ has degree at most $2g$.

Notice that if $i > p(2k+1)/2+pm$, then $A_i(x)$ is divisible by $p^m$,
and by Lemma~\ref{lem:est}, the reduction of $A_i(x)y^{-2i-1}\,dx$
will be divisible by $p^{m - \lfloor \log_p (2m+1)\rfloor}$.
Therefore the reduction of
$(A(x)y^{-2k-1}\,dx)^\sigma$
is determined by the $A_i$ with $i \leq
N_1 + \log_p (2N_1)$.

\section{An algorithm for computing Frobenius}

Using the results of the previous section, we now describe an algorithm
for computing the characteristic polynomial of Frobenius on a
hyperelliptic curve $C$ of genus $g$ over $\FF_q$ with $q=p^n$.
We maintain the notation of the previous section.

\subsection*{Step 1: Compute Frobenius on $y$}
Compute a sequence of polynomials $A_0(x), A_1(x), A_{pN-1}(x)$ over
$W(\FF_q)/(p^{N_1})$, each of degree at most $2g$, such that
\begin{align*}
\frac{1}{y^\sigma} &= y^{-p} \left(1 + \frac{Q(x)^\sigma - Q(x)^p}{y^{2p}}\right)^{-1/2} \\
&\equiv y^{-p} \sum_{i=0}^{pN-1} \frac{A_i(x)}{y^{2i}}
\end{align*}
as a power series in $y^{-2}$ over $W(\FF_q)/(p^{N_1})$ modulo $y^{-2pN}$ 
using a Newton iteration. Specifically, recall that for
$s \in 1 + tK[[t]]$, to compute $s^{-1/2}$ we may
set $x_0 = 1$ and
\[
x_{i+1} \equiv \frac{3}{2} x_i - \frac{1}{2} sx_i^3 \pmod{t^{2^{i+1}}};
\]
then $x_i \equiv s^{-1/2} \pmod{t^{2^i}}$.
The dominant operation in this iteration is the cubing, which can be done
in asymptotically optimal time by, for instance,
packing $x_i$ into an integer and applying
the Sch\"onhage-Strassen algorithm for fast integer multiplication.

\subsection*{Step 2: Compute Frobenius on differentials}
For $i=0, \dots, 2g-1$, compute the reduction of $(x^i\,dx/y)^\sigma$ as follows.
Using the computation of $1/y^\sigma$ carried out in the first step, write
\[
\left(\frac{x^i\,dx}{y}\right)^\sigma
= \frac{px^{pi+(p-1)}\,dx}{y^\sigma}
= \frac{G(x)\,dx}{y} +
\sum_{j=1}^{pN}
\frac{F_i(x)\,dx}{y^{2j+1}} + O(y^{-2pN-3}),
\]
where $\deg F_i \leq 2g-1$; for notational convenience, set $F_0(x) = 0$.
Then compute $S_k(x)$ for $k=2pN, 2pN-1, \dots, 1$ as follows.
Let $S_{2pN}(x) = F_{2pN}(x)$. Given $S_{k+1}(x)$, find polynomials
$A_{k+1}(x)$ and $B_{k+1}(x)$ such that $A_{k+1}Q + B_{k+1}Q' = S_{k+1}$.
Then set $S_k(x) = F_{k} + A_{k+1} + 2B_{k+1}/(2k-1)$. 
By
the reduction argument from the previous section, $(x^i\,dx/y)^\sigma$
is cohomologous to $(S_0(x)+G(x))\,dx/y$.

Note that the above computation cannot be performed in
$W(\FF_q)/(p^N)$ as written, because of the division by $2k-1$.
To remedy this, interpret $2B_{k+1}/(2k-1)$ to mean any
polynomial over $W(\FF_q)/(p^N)$ which, when multiplied by $2k-1$,
equals $2B_{k+1}$. Lemma~\ref{lem:est} implies both that any discrepancy
introduced in $S_k$ in case $2k-1$ is divisible by $p$
has no effect on $S_1$ modulo $p^{N_1}$, and that $B_{k+1}/(2k-1)$ always
has integral coefficients.

By construction, $S_0$ has degree at most $2g$, but $G$ can have
degree up to $2pg-1$, so we must reduce $G(x)\,dx/y$ in cohomology
as well. For $j = \deg G - 2g+1$, set $G_j(x) = G(x)$; for 
$k=j,j-1,\dots,1$, let $G_{k-1}$ be the remainder of $G_k(x)$
modulo $x^{k-1} Q'(x) + 2 (k-1) x^{k-2} Q(x)$. (When the latter
has leading coefficient divisible by $p$, we may again fill the high 
$p$-adic digits arbitrarily without affecting the final result of the
computation, by Lemma~3.)
Then
$G(x)\,dx/y$ is cohomologous to $G_0(x)\,dx/y$, and so
$(x^i\,dx/y)^\sigma$ is cohomologous to $(S_1+G_0)\,dx/y$.

\subsection*{Step 3: Compute characteristic polynomial}

From the previous step, we may extract the matrix $M$ 
through which the $p$-power Frobenius acts on a basis of cohomology
over $W(\FF_q)/(p^N)$.
Compute $M' = M M^\sigma M^{\sigma^2} \cdots M^{\sigma^{n-1}}$,
determine the characteristic polynomial of $M'$, and recover
the characteristic polynomial of Frobenius from the first $g$ coefficients
modulo $p^N$.

In case one wants only the Newton polygon of Frobenius and not its
full characteristic polynomial, some savings may be possible in this step.
On one hand, although the Newton polygon can be computed directly from the
$p$-Frobenius, it is not given by the characteristic polynomial of
the matrix $M$ in general. On the other hand, this does work in case $M=DA$
with $D$ diagonal and $A$ congruent to the identity matrix modulo $p$, and
in some cases it may be easy to select a basis for which this holds.

\section{Resource analysis}

We now analyze the space and time requirements of the algorithm
for a curve of genus $g$ over $\FF_{p^n}$ (keeping $p$ fixed).
Before proceeding through the individual steps, we make some general
observations about the implementation of low-level operations that
permeate the discussion.

All ring operations in the algorithm
take place in the degree $n$ unramified extension
of $\ZZ_p/(p^N)$, 
and each element of this ring requires $O(gn^2)$ storage space.
Using fast integer multiplication as noted above,
individual multiplications and divisions
in the ring can be accomplished in time $O(g^{1+\epsilon} n^{2+\epsilon})$.

Applying any power $\tau = \sigma^k$ of the ring automorphism $\sigma$
can be accomplished in time
$O(g^{1+\epsilon} n^{3+\epsilon})$ as follows. Suppose the base ring
is represented as $\ZZ_p/(p^N)[\alpha]$ where $P(\alpha) = 0$. Compute
an element of the residue field congruent to $\alpha^{p^k}$ mod $p$
by repeated squarings. Then use Newton's iteration to compute
$\alpha^{\tau}$ from this. Now to compute $G(\alpha)^{\tau}$,
for $G$ a polynomial over $\ZZ_p/(p^N)$, evaluate $G$ at $\alpha^{\tau}$
using Horner's method, or
(better in practice) the Paterson-Stockmeyer algorithm \cite{ps},
using $O(n)$ multiplications in $\ZZ_p(p^N)$.

In Step~1, we compute $O(gn)$ terms of $1/y^\sigma$; each term
consists of a polynomial in $x$ of degree at most $2g-1$, which
requires $O(g^2n^2)$ space to  store. Thus the entire expression
requires space $O(g^3n^3)$ and time
$O(g^{3+\epsilon} n^{3+\epsilon})$ to compute.

In Step~2, the dominant step in each reduction is writing a polynomial $T$
of degree at most $2g-1$ as $AQ + BQ'$. This can be done by precomputing
polynomials
$R$ and $S$ of degrees $2g-1$ and $2g$, respectively, such that $RQ + SQ' = 1$,
then computing $A$ as the reduction of $TQ$ modulo $Q'$ and $B$ as the
reduction of $SQ'$ modulo $Q$. Since the polynomials in question require
space $O(g^{2+\epsilon} n^{2+\epsilon})$ each, this extended GCD operation
can be performed in time $O(g^{2+\epsilon} n^{2+\epsilon})$;
see \cite{moenck}.
The reduction step is performed $O(gn)$ times for each of $2g$ forms,
for a total of $O(g^{4+\epsilon} n^{3+\epsilon})$ time.

In Step~3, we begin with a $2g \times 2g$ matrix $M$ each of whose entries
has size $O(gn^2)$, and must compute
$M' = M M^\sigma \cdots M^{\sigma^{n-1}}$ by
repeated squaring. Specifically, we can compute $M_1 = MM^\sigma$,
$M_2 = M_1M_1^{\sigma^2}$, $M_3 = M_2M_2^{\sigma^4}$ and so on,
then combine these as in the usual repeated squaring method for
exponentation to compute $M'$.
This process requires $O(\log n)$ multiplications of $2g \times 2g$
matrices and $O(g^2 \log n)$ applications of powers of $\sigma$
(specifically, of powers of the form $\sigma^m$ for $m$ a power of 2).
The former requires $O(g^3 \log n)$ ring operations,
at a cost of $O(g^{4+\epsilon} n^{2+\epsilon})$ time;
the latter requires $O(g^{3+\epsilon} n^{3+\epsilon})$ time.

We then must compute the characteristic polynomial of $M'$.
This can be accomplished in $O(g^3)$ ring operations, e.g.,
by computing $v, Mv, M^2v, \dots$ until these fail to be linearly
independent, then inverting a matrix to obtain a factor of the characteristic
polynomial, and repeating as needed.
This translates into a time cost of $O(g^{4+\epsilon} n^{2+\epsilon})$.

Overall, the dominant factors are $g^{4+\epsilon}$ and $n^{3+\epsilon}$.
Note, however, that one factor of $g$ can be saved in a parallel computation
in Step~2,
by computing the Frobenius on each basis vector simultaneously.
On the other hand, the factor $g^{4+\epsilon}$ remains as a bottleneck
in Step~3 and does not appear to be readily mollifiable by parallelism.
Likewise, a parallel approach does not appear to mollify the 
factor of $n^{3 + \epsilon}$ appearing throughout the
analysis.

\section*{Acknowledgments}
Thanks to Dan Bernstein for suggesting Moenck's extended GCD
algorithm, to Takakazu Satoh for suggesting the Paterson-Stockmeyer
method, and to Johan de~Jong,
Joe Wetherell and Hui June Zhu for helpful discussions.
The author was supported by an NSF Postdoctoral Fellowship.


\begin{thebibliography}{9}

\bibitem{berthelot}
P. Berthelot, Finitude et puret\'e cohomologique en cohomologie rigide
(with an appendix in English by Aise Johan de Jong),
\textit{Invent. Math.} \textbf{128} (1997), 329--377. 

\bibitem{cohen}
H. Cohen, \textit{A course in computational algebraic number theory},
Graduate Texts in Mathematics \textbf{138}, Springer-Verlag, 1993.

\bibitem{fgh}
M. Fouquet, P. Gaudry, and R. Harley, An extension of Satoh's algorithm
and its implementation, \textit{J. Ramanujan Math. Soc.} \textbf{15} (2000),
281--318.

\bibitem{hartshorne}
R. Hartshorne, \textit{Algebraic Geometry}, Graduate Texts in Mathematics
52, Springer-Verlag (New York), 1977.

\bibitem{lw}
A.G.B. Lauder and D. Wan, Counting points on varieties over finite
fields of small characteristic, preprint (available at \verb+www.math.uci.edu/~dwan+).

\bibitem{moenck}
R.T. Moenck, Fast computation of GCDs, in
\textit{Fifth Annual ACM Symposium on Theory of Computing
(Austin, Tex., 1973)}, Assoc. Comput. Mach. (New York), 1973, pp. 142--151.

\bibitem{mw}
P. Monsky and G. Washnitzer, Formal cohomology. I, \textit{Ann. of Math. (2)}
\textbf{88} (1968), 181--217.

\bibitem{m2}
P. Monsky, Formal cohomology. II. The cohomology sequence of a pair,
 \textit{Ann. of Math. (2)} \textbf{88}
(1968), 218--238.

\bibitem{m3}
P. Monsky, Formal cohomology. III. Fixed point theorems,
\textit{Ann. of Math. (2)} \textbf{93} (1971), 315--343.

\bibitem{ps}
M.S. Paterson and L.J. Stockmeyer,
On the number of nonscalar multiplications necessary to evaluate polynomials,
\textit{SIAM J. Comput.} \textbf{2} (1973), 60--66.

\bibitem{pila}
J. Pila, Frobenius maps of abelian varieties
and finding roots of unity in finite fields,
\textit{Math. Comp.} \textbf{55} (1990), 745--763.

\bibitem{satoh}
T. Satoh, The canonical lift of an ordinary elliptic curve over a finite field
and its point counting, \textit{J. Ramanujan Math. Soc.} \textbf{15} (2000),
247--270.


\bibitem{schoof}
R. Schoof, Elliptic curves over finite fields and the computation of
square roots mod $p$, \textit{Math. Comp.} \textbf{44} (1985), 483--494.

\bibitem{vdp}
M. van~der~Put, The cohomology of Monsky and Washnitzer, in
Introductions aux cohomologies
$p$-adiques (Luminy, 1984), \textit{M\'em. Soc. Math. France} \textbf{23}
(1986), 33--60.

\end{thebibliography}
\end{document}